\documentclass[10pt]{article}

\usepackage[utf8]{inputenc}
\usepackage[T1]{fontenc}
\usepackage{mathptmx}
\usepackage{mathtools}
\usepackage{stmaryrd}
\usepackage{amsmath}
\usepackage{amssymb}
\usepackage{amsthm}
\usepackage{tikz-cd}
\usepackage{geometry}
\usepackage{pgfplots}
 
\pgfplotsset{compat = newest}

\newtheorem{thm}{Theorem}[section]
\newtheorem{prop}[thm]{Proposition}

\newtheorem{lmm}[thm]{Lemma}

\newtheorem{conj}[thm]{Conjecture}

\theoremstyle{definition}

\newtheorem*{rmk}{Remark}

\newenvironment{prf}{{\noindent\it Proof}\quad}{\hfill $\square$\par}

\title{On some birational invariants of hyper-Kähler manifolds}
\author{Chenyu Bai}
\date{}

\begin{document}
\maketitle

\begin{abstract}
    We study in this article three birational invariants of projective hyper-Kähler manifolds: the degree of irrationality, the fibering gonality and the fibering genus. We first improve the lower bound in a recent result of Voisin saying that the fibering genus of a Mumford--Tate very general projective hyper-Kähler manifold is bounded from below by a constant depending on its dimension and the second Betti number. We also study the relations between these birational invariants for projective $K3$ surfaces of Picard number $1$ and study the asymptotic behaviors of their degree of irrationality and fibering genus.
\end{abstract}

\section{Introduction}
In the paper~\cite{Irrationality}, a birational invariant, the degree of irrationality, of projective manifolds is studied, in particular for hypersurfaces. The degree of irrationality $\mathrm{Irr}(X)$ of a projective manifold $X$ is defined to be the minimal degree of all possible dominant rational maps from $X$ to the projective space $\mathbb P^n$, with $n={\dim X}$. The degree of irrationality is a natural generalization of gonality of curves, and it measures to which extent $X$ is far from being rational, in the sense that $X$ is rational if and only if $\mathrm{Irr}(X)=1$. 

The degree of irrationality is now relatively well-understood for very general hypersurfaces in projective spaces~\cite{Irrationality, ThesisOfStapleton, Yang} and for very general abelian varieties~\cite{IrrAVVoisin, ThesisOfMartin}, while for $K3$ surfaces, it is still mysterious. In particular, the following conjecture proposed in~\cite{Irrationality} remains open and is a main motivation of the second part of this paper. 

\begin{conj}[\cite{Irrationality}]\label{Conjecture}
Let $\{(S_d, L_d)\}_{d\in \mathbb N}$ be very general polarized $K3$ surfaces such that $L_d^2=2d-2$. Then
\[\limsup_{d\to \infty}\mathrm{Irr}(S_d)=+\infty.
\]
\end{conj}

In an attempt to understand and estimate the degree of irrationality, many other birational invariants are proposed among which are the fibering gonality and fibering genus introduced in~\cite{FibgenFibgon}. As is defined in \emph{loc. cit.}, for a projective manifold $X$, the fibering gonality $\mathrm{Fibgon}(X)$ (resp. fibering genus $\mathrm{Fibgen}(X)$) of $X$ is the minimal gonality (resp. geometric genus) of the general fiber of a fibration $X\dashrightarrow B$ into curves.

In this article, we study these three birational invariants, namely the \emph{degree of irrationality}, the \emph{fibering gonality} and the \emph{fibering genus}, for projective hyper-Kähler manifolds. A projective hyper-Kähler manifold $X$ is a simply connected complex projective manifold such that $H^0(X,\Omega_X^2)=\mathbb C\sigma_X$, where $\sigma_X$ is a nowhere degenerate holomorphic $2$-form on $X$. A projective $K3$ surface is a $2$-dimensional projective hyper-Kähler manifold.

Let $X$ be a projective hyper-Kähler manifold of dimension $2n$. Let $q_X$ be its Beauville--Bogomolov--Fujiki form. 
Let $b_{2,tr}(X)$ be the dimension of $H^2(X,\mathbb Q)_{tr}$, the transcendental part of $H^2(X,\mathbb Q)$. In~\cite[Theorem 1.5]{FibgenFibgon}, Voisin proves the following result on the fibering genus of Mumford--Tate very general projective hyper-Kähler manifolds. 

\begin{thm}[Voisin \cite{FibgenFibgon}]\label{TheoremFibGenVoisin}
 Assume that the Mumford--Tate group of the Hodge structure on $H^2
(X, \mathbb Q)_{tr}$ is maximal. Assume that $n \geq 3$ and $b_2(X)_{tr} \geq
5$. 
Then if $X$ admits a fibration $\phi: X\dashrightarrow B$, with $\dim B = 2n - 1$, the general fiber of $\phi$ has
genus $g \geq \mathrm{Inf}(n + 2, 2^{\lfloor\frac{ b_{2,tr}(X)-3}{2}\rfloor})$. In other words,
$\mathrm{Fibgen}(X) \geq \mathrm{Inf}(n + 2, 2^{\lfloor\frac{ b_{2,tr}(X)-3}{2}\rfloor})$.
\end{thm}

In this theorem, the assumption on the maximality of the Mumford--Tate group is equivalent to saying that it is the special orthogonal group of $(H^2(X, \mathbb Q)_{tr}, q_X)$. This assumption is satisfied by the very general fiber of a complete lattice polarised family of hyper-Kähler manifolds.

Our first result is a slight improvement of this theorem, giving a sharper lower bound of the fibering genus of Mumford--Tate very general projective hyper-Kähler manifolds. This improvement is obtained by a more refined discussion on an orthogonal relation (see Proposition~\ref{LemmeDeVoisin} (iii) for details) observed in~\cite{FibgenFibgon}.

\begin{thm}\label{TheoremFibGen}
 Assume that the Mumford--Tate group of the Hodge structure of $H^2(X,\mathbb Q)_{tr}$ is maximal. Assume that $b_{2,tr}(X)\geq 5$. Then 
\[\mathrm{Fibgen}(X)\geq \min\left\{ n+\left\lceil\frac{-1+\sqrt{8n-7}}{2}\right\rceil, 2^{\lfloor\frac{ b_{2,tr}(X)-3}{2}\rfloor}
\right\}.
\]
\end{thm}

Our second result is about the relations between the degree of irrationality, the fibering gonality and the fibering genus of a projective K3 surface with Picard number $1$.
\begin{thm}\label{ThmFibgonDeK3}
Let $S$ be a projective K3 surface whose Picard number is $1$. Then one of the following two cases holds:\\
(a) $\mathrm{Irr}(S)=\mathrm{Fibgon}(S)$;\\
(b) $\mathrm{Fibgen}(S)^2\leq \mathrm{Fibgon}(S)^{21}$.
\end{thm}
In particular, if the fibering gonality of $S$ is very small compared to the fibering genus, then it is equal to the degree of irrationality of $S$. The fibering genus of a projective K3 surface of Picard number $1$ is studied in~\cite{EinLazarsfeld}. We will rely on~\cite[Theorem B]{EinLazarsfeld} (stated below as Theorem~\ref{EinLazarsfeld} for completeness).
\begin{thm}[Ein--Lazarsfeld~\cite{EinLazarsfeld}]\label{EinLazarsfeld}
Let $S$ be a projective $K3$ surface whose Picard group is generated by a line bundle of degree $2d-2$. Then the fibering genus of $S$ is of order $\sqrt{d}$. More precisely,
\begin{equation}\label{BoundsForFibgenofS}
    \sqrt{\frac{d}{2}}\leq \mathrm{Fibgen}(S)\leq 2\sqrt{2d}.
\end{equation}
\end{thm}
 Using Theorem~\ref{ThmFibgonDeK3} and~\cite[Theorem B]{EinLazarsfeld}, we can compare the asymptotic behaviors of the degree of irrationality and the fibering gonality of very general polarized $K3$ surfaces. 
\begin{thm}\label{TendsToInfty}
Let $\{S_d\}_{d\in\mathbb N}$ be projective K3 surfaces such that the Picard group of $S_d$ is generated by a line bundle with self intersection number $2d-2$. Then 
\[\limsup_{d\to\infty} \mathrm{Irr}(S_d)=+\infty \iff \limsup_{d\to \infty} \mathrm{Fibgon}(S_d)=+\infty.\]
\end{thm}
\begin{proof}[Proof of Theorem~\ref{TendsToInfty} assuming Theorem~\ref{ThmFibgonDeK3}]
The implication $\Leftarrow$ is clear since $\mathrm{Fibgon}(S_d)\leq \mathrm{Irr}(S_d)$ for any $d$. Now let us prove the implication $\Rightarrow$. Theorem~\ref{EinLazarsfeld} (\cite[Theorem B]{EinLazarsfeld}) shows that $\lim_{d}\mathrm{Fibgen}(S_d)=+\infty$. If $\limsup_{d} \mathrm{Fibgon}(S_d)\neq+\infty$, then there exist constants $C$ and $D$ such that for every $d>D$ we have $\mathrm{Fibgon}(S_d)^{21/2}<C<\mathrm{Fibgen}(S_d)$. Hence for $d>D$, we must have $\mathrm{Irr}(S_d)=\mathrm{Fibgon}(S_d)<C^{2/21}$ by Theorem~\ref{ThmFibgonDeK3}, which contradicts the assumption that $\limsup_{d} \mathrm{Irr}(S_d)=+\infty$.
\end{proof}

Conjecture~\ref{Conjecture} predicts that $\limsup_{d} \mathrm{Irr}(S_d)=+\infty$, so it should be expected that $\limsup_{d} \mathrm{Fibgon}(S_d)=+\infty$ as well.

In the course of the proof of Theorem~\ref{ThmFibgonDeK3}, we obtain the following inequalities about the discriminant of Picard lattices of projective hyper-Kähler manifolds. Although we only use this proposition in the case of $K3$ surfaces and we only use one side of the inequalities, the following general form is of independent interest. Here $\mathrm{disc}(\mathrm{Pic}(X))$ denotes the discriminant of the Picard lattice of $X$ with respect to the Beauville--Bogolov--Fujiki form and similarly for $\mathrm{disc}(\mathrm{Pic}(X'))$.
\begin{prop}\label{InequalityOfLatticesHK}
Let $X$ and $X'$ be deformation equivalent projective hyper-Kähler manifolds of dimension $2n$. Let $\phi: X\dashrightarrow X'$ be a dominant rational map. Then
\[(\deg\phi)^{(\frac{1}{n}-2)b_{2,tr}(X)}\leq \left|\frac{\mathrm{disc}(\mathrm{Pic}(X))}{\mathrm{disc}(\mathrm{Pic}(X'))}\right|\leq (\deg\phi)^{\frac{b_{2,tr}(X)}{n}}.
\]
\end{prop}
If  in the previous proposition $X$ and $X'$ are $K3$ surfaces, we can get a slightly better lower bound.
\begin{prop}\label{InequalityOfLatticesK3}
Let $S$ and $S'$ be projective $K3$ surfaces and let $\lambda(S)=\min\{\rho(S), b_{2,tr}(S)\}$ where $\rho(S)$ is the Picard number of $S$. Let $\phi: S\dashrightarrow S'$ be a dominant rational map. Then
\[(\deg\phi)^{-\lambda(S)}\leq \left|\frac{\mathrm{disc}(\mathrm{Pic}(S))}{\mathrm{disc}(\mathrm{Pic}(S'))}\right|\leq (\deg\phi)^{b_{2,tr}(S)}.
\]
\end{prop}

The organization of the article is as follows. 
In Section~\ref{FibgenOfHK}, we prove Theorem~\ref{TheoremFibGen}. In Section~\ref{K3}, we prove Theorem~\ref{ThmFibgonDeK3}. In Section~\ref{Inequalities}, we prove Propositions~\ref{InequalityOfLatticesHK} and \ref{InequalityOfLatticesK3}.

\subsection*{Acknowledgment}
It is a pleasure to thank my PhD advisor Claire Voisin for sharing her ideas with me. This work is initiated and motivated by her interesting paper~\cite{FibgenFibgon}. I would like to thank her specially for fruitful and patient discussions on her paper. I would also like to thank Fabrizio Anella and Mauro Varesco for helpful discussions on curve theory, and Pietro Beri and Yicheng Zhou for interesting discussions on lattices. I am also grateful to Olivier Martin, David Stapleton, Robert Lazarsfeld and Angel David Rios Ortiz for encouraging me to work on measures of irrationality. 

The conference \emph{Algebraic Geometry} in honor of Claire Voisin's 60th birthday and the conference \emph{Géométrie Algébrique en Liberté XXIX} gave me lots of occasions to meet people and to discuss this work. I would like to thank the organizers of these two beautiful conferences. This work was done during the preparation of my PhD thesis. I would like to thank the \emph{Institut de Mathématiques de Jussieu - Paris Rive Gauche} for marvelous research environment. The thesis is supported by the ERC Synergy Grant HyperK (Grant agreement No. 854361).

\section{Fibering genus of very general hyper-Kähler manifolds}\label{FibgenOfHK}
 Let $X$ be a projective hyper-Kähler manifold of dimension $2n$. Let $f: X\dashrightarrow B$ be a fibration into curves and let $\tau: \tilde X\to X$ and $\tilde f: \tilde X\to B$ be a resolution of indeterminacy points. Let $\tilde X_b$ be a smooth fiber of $\tilde f$ over a general point $b\in B$. Contracting with $\tau^*\sigma_{X}$ gives a morphism of vector bundles $N_{\tilde X_b/\tilde X}\to \Omega_{\tilde X_b}$. Therefore, we have a morphism
\[\sigma_b: T_{B,b}\to H^0(\tilde X_b,
\Omega_{\tilde X_b})=H^0(\tilde X_b, K_{\tilde X_b})\] defined as the composition of the natural morphisms $T_{B, b}\to H^0(\tilde X_b, N_{\tilde X_b/\tilde X})$ and $\lrcorner \tau^*\sigma_X: H^0(\tilde X_b, N_{\tilde X_b/\tilde X})\to H^0(\tilde X_b,\Omega_{\tilde X_b})$. Let $\rho: T_{B,b}\to H^1(\tilde X_b,T_{\tilde X_b})$ be the Kodaira--Spencer map. In~\cite{FibgenFibgon}, Voisin proves the following

\begin{prop}[Voisin~\cite{FibgenFibgon}]\label{LemmeDeVoisin}
 (i) If $n\geq 2$, then $\tilde X_b$ is \emph{not} hyperelliptic.\\
 (ii) The rank of $\sigma_b$ is $\geq n$.\\
 (iii) Let $I_b$ be image of $\mathrm{Im}\sigma_b\otimes H^0(\tilde X_b,K_{\tilde X_b})$ under the multiplication map
 \[\mu: H^0(\tilde X_b,K_{\tilde X_b})\otimes H^0(\tilde X_b,K_{\tilde X_b})\to H^0(\tilde X_b,K_{\tilde X_b}^{\otimes 2}).\] Then $\rho(\ker\sigma_b)\subset H^1({\tilde X_b}, T_{\tilde X_b})$ is orthogonal to $I_b$ via the Serre pairing
 \[ H^1({\tilde X_b}, T_{\tilde X_b})\otimes H^0({\tilde X_b},K_{\tilde X_b}^{\otimes 2})\to H^1({\tilde X_b}, K_{\tilde X_b})\cong \mathbb C.\]
  (iv) Assume that that $b_{2,tr}(X)\geq 5$ and that the Mumford--Tate group of $H^2(X,\mathbb Q)_{tr}$ is maximal. If $\rho: T_{B,b}\to H^1({\tilde X_b},T_{\tilde X_b})$ is \emph{not} injective, then $g({\tilde X_b})\geq 2^{\lfloor\frac{ b_{2,tr}(X)-3}{2}\rfloor}$.
 \end{prop}
 The proofs of (i), (ii) and (iii) are explicitly written in~\cite{FibgenFibgon}. Although (iv) is not stated in~\cite{FibgenFibgon} with this generality, it is essentially proved there (see the proof of Lemma 2.13 and Lemma 2.15 in \textit{loc. cit.}).

\subsection{Proof of Theorem~\ref{TheoremFibGen}}
For $n=1$,  the inequality is obvious, since  $K3$ surfaces cannot be fibered into rational curves. From now on, let us assume $n\geq 2$ to be able to assume $\tilde X_b$ is not hyperelliptic (see Proposition~\ref{LemmeDeVoisin} (i)). Let $k$ be the corank of $\sigma_b$, i.e., $\mathrm{rank}\sigma_b=g-k$, where $g$ is the genus of the curve $\tilde X_b$. We now prove 
\begin{lmm}\label{LinearAlgebraLemma}
The notation as above, we have the following inequality $$\dim\ker\rho \geq 2n-1-(g-k)-\frac{k(k+1)}{2}.$$
\end{lmm}
\begin{proof}
We use the notation as in Proposition~\ref{LemmeDeVoisin} (iii). By Serre duality, the orthogonality result in Proposition~\ref{LemmeDeVoisin} (iii) implies that $\dim \rho(\ker\sigma_b)+\dim I_b\leq \dim  H^0({\tilde X_b},K_{\tilde X_b}^{\otimes 2} )$, that is,
\begin{equation}\label{Inequality1inFibgen}
    \mathrm{codim}(I_b\subset H^0({\tilde X_b},K_{\tilde X_b}^{\otimes 2} ))\leq \dim \rho(\ker\sigma_b).
\end{equation}
The multiplication map $\mu: H^0(\tilde X_b,K_{\tilde X_b})\otimes H^0(\tilde X_b,K_{\tilde X_b})\to H^0(\tilde X_b,K_{\tilde X_b}^{\otimes 2})$ factors through the symmetric product $\mathrm{Sym}^2H^0(\tilde X_b,K_{\tilde X_b})$, that is, we have a commutative diagram
\[\begin{tikzcd}
H^0(\tilde X_b,K_{\tilde X_b})\otimes H^0(\tilde X_b,K_{\tilde X_b}) \arrow[rr, "\mu"]\arrow[rd, "\mathrm{pr}"]&& H^0(\tilde X_b,K_{\tilde X_b}^{\otimes 2})\\
& \mathrm{Sym}^2H^0(\tilde X_b,K_{\tilde X_b})\arrow[ru, "\mu'"]&
\end{tikzcd},
\]
where $\mathrm{pr}:H^0(\tilde X_b,K_{\tilde X_b})\otimes H^0(\tilde X_b,K_{\tilde X_b})\to \mathrm{Sym}^2H^0(\tilde X_b,K_{\tilde X_b})$ is the canonical symmetrization map. Let \mbox{$\mathrm{Im}\,\sigma_b\cdot H^0(\tilde X_b,K_{\tilde X_b})$} denote the image of $\mathrm{Im}\,\sigma_b\otimes H^0(\tilde X_b,K_{\tilde X_b})$ under $\mathrm{pr}$. Then we have 
\begin{equation}\label{Equality1inFibgen}
\mathrm{codim}(\mathrm{Im}\,\sigma_b\cdot H^0(\tilde X_b,K_{\tilde X_b}))\subset \mathrm{Sym}^2H^0(\tilde X_b,K_{\tilde X_b}))
=\frac{k(k+1)}{2}.
\end{equation}
On the other hand, by Max Noether theorem (see~\cite[Chapter III, §2]{Arbarello} or~\cite{Voisin}), the multiplication map $\mu': \mathrm{Sym}^2H^0(\tilde X_b,K_{\tilde X_b})\to H^0(\tilde X_b,K_{\tilde X_b}^{\otimes 2})$ is surjective, since $\tilde X_b$ is not hyperelliptic. Taking into account the fact that $I_b$ is the image of $\mathrm{Im}\,\sigma_b\cdot H^0(\tilde X_b,K_{\tilde X_b})$ under the map $\mu'$, we have the following inequality
\begin{equation}\label{Inequality2inFibgen}
\mathrm{codim}(I_b\subset H^0(\tilde X_b,K_{\tilde X_b}^{\otimes 2}) ) \leq \mathrm{codim}(\mathrm{Im}\,\sigma_b\cdot H^0(\tilde X_b,K_{\tilde X_b})\subset \mathrm{Sym}^2H^0(\tilde X_b,K_{\tilde X_b})).
\end{equation}
Combining (\ref{Inequality1inFibgen}), (\ref{Equality1inFibgen}) and (\ref{Inequality2inFibgen}), we get $\dim\rho(\ker\sigma_b)\leq \frac{k(k+1)}{2}$, from which we deduce that $$\dim\ker\rho\geq \dim\ker\sigma_b-\frac{k(k+1)}{2}=2n-1-(g-k)-\frac{k(k+1)}{2}.$$
\end{proof}  

\begin{proof}[Proof of Theorem~\ref{TheoremFibGen}] Assuming  that $g(\tilde X_b)< 2^{\lfloor\frac{ b_{2,tr}(X)-3}{2}\rfloor}$, we have to prove $g\geq n+ \lceil\frac{-1+\sqrt{8n-7}}{2}\rceil$. By Proposition~\ref{LemmeDeVoisin} (ii), (iv) and Lemma~\ref{LinearAlgebraLemma}, we have the following constraints on $g$ and $k$:
\[\left\{\begin{array}{lr}
     k\geq 0\\
     & \cr 
     g-k- n\geq 0  \\
      & \cr 
     2n-1-(g-k) -\frac{k(k+1)}{2} \leq 0.
\end{array}\right.
\]
In order to find the minimal possible value of $g$ under these constraints, we make the following discussion according to the values of $k$. 
\begin{itemize}
    \item When $0\leq k\leq \frac{-1+\sqrt{8n-7}}{2}$, we have $k-\frac{k(k+1)}{2}+2n-1\geq n+k$. Hence, the minimal possible value of $g$ in this domain is the minimum of $k-\frac{k(k+1)}{2}+2n-1$ with $0\leq k\leq  \frac{-1+\sqrt{8n-7}}{2}$, which is $n+ \frac{-1+\sqrt{8n-7}}{2}$.
    \item When $k\geq  \frac{-1+\sqrt{8n-7}}{2}$, we have $k-\frac{k(k+1)}{2}+2n-1\leq n+k$. Hence, the minimal possible value of $g$ in this domain is the minimum of $n+k$ with  $k\geq  \frac{-1+\sqrt{8n-7}}{2}$, which is $n+ \frac{-1+\sqrt{8n-7}}{2}$.
\end{itemize}

Since $g$ and $k$ are integers, we find $g\geq n+ \lceil\frac{-1+\sqrt{8n-7}}{2}\rceil$, as desired.

\end{proof}

\begin{rmk}
Our proof relies on the inequality (\ref{Inequality2inFibgen}) which only uses the surjectivity of the multication map $\mu'$. With more information on the geometry of the canonical embedding, and in particular, on the gonality of the fibers, we could get a better estimate in Theorem~\ref{TheoremFibGen}.

\end{rmk}

\section{Relations between birational invariants of $K3$ surfaces}\label{K3}
In this section, we are going to prove Theorem~\ref{ThmFibgonDeK3} that relates the three birational invariants, namely, the degree of irrationality, the fibering gonality and the fibering genus, of projective $K3$ surfaces of Picard number $1$.

\subsection{A factorization}\label{ConstructionofPhi}
Let $S$ be a smooth projective surface and let $f: S\dashrightarrow B$ be a fibration into curves over a smooth base $B$. After a resolution of inderminacies of $f$ and replacing $S$ by another birational model, we may assume $f: S\to B$ is a morphism. Let $d$ be the gonality of the general fiber of $f$. The general fiber $C$ admits a degree $d$ morphism from $C$ to $\mathbb P^1$. Standard argument shows that we can spread this morphism into a family up to a generically finite base change.
\begin{lmm}\label{StandardArgumentsGonalityK3}
There is a generically finite morphism $\pi: B'\to B$ and a degree $d$ dominant rational map $\psi: S\times_BB'\dashrightarrow \mathbb P^1\times B'$ over $B'$.
\end{lmm}
\begin{prf}
Let $B_0$ be the smooth locus of $f: S\to B$ and let $f_0: S_0\to B_0$ be the restriction of $f$ on smooth locus. Let $p: \mathrm{Pic}^d(S_0/B_0)\to B_0$ be the degree $d$ relative Picard variety of $f_0: S_0\to B_0$. By the assumption on the general fiber of $f: S\to B$, the Brill-Noether locus in $\mathrm{Pic}^d(S_0/B_0)$ of linear systems of degree $d$ and dimension $1$ is dominant over $B_0$ via the map $p: \mathrm{Pic}^d(S_0/B_0)\to B_0$.  Let $B_0'$ be a general reduced irreducible subscheme of $\mathrm{Pic}^d(S_0/B_0)$ that is dominant and generically finite over $B_0$ by $p$. Let us take $B'$ to be a completion of $B_0'$. Then by construction, the universal line bundle restricted to $S_0\times_BB_0'$ gives a dominant rational map $\psi: S\times_BB'\dashrightarrow \mathbb P_{B'}^1$ of degree $d$, as desired.
\end{prf}

It is natural to ask if $\psi: S\times_BB'\dashrightarrow \mathbb P^1\times B'$ over $B'$ descends to a rational map $\psi_B: S\dashrightarrow \mathbb P^1\times B$ over $B$. A moment of thinking will convince us that we are asking too much, because $\pi: B'\to B$ is in general not a Galois cover. We make the following construction. Let $n$ be the degree of the morphism $\pi: B'\to B$. Consider the $n$-th self fibred product of $B'$ over $B$: $B'\times_B\ldots\times_BB'$. Define $B''$ to be the closure in $B'\times_B\ldots\times_BB'$ of the set 
\[\{(x_1, \ldots, x_n)\in B'\times_B\ldots\times_BB': x_1,\ldots, x_n 
\textrm{ are distinct }\}.
\]
Then $\pi': B''\to B$ is of degree $n!$ and the symmetric group $\mathfrak S_n$ permuting the components of $B'\times_B\ldots\times_BB'$ acts on an open dense subset of $B''$. The rational map $\psi: S\times_BB'\dashrightarrow \mathbb P^1\times B'$ over $B'$ given in Lemma~\ref{StandardArgumentsGonalityK3} can be extended to a rational map 
\begin{equation}\label{EqPsi}
\psi': S\times_BB''\dashrightarrow (\mathbb P^1)^n\times B''
\end{equation}
over $B''$ in a natural way: let $x\in S$ and let $y=(y_1,\ldots, y_n)\in B''$ be general points, we define $\psi'(x,y)=(\psi(x,y_1),\ldots, \psi(x,y_n), y)$. Moreover, the symmetric group $\mathfrak S_n$ acts canonically on both sides of (\ref{EqPsi}) in the following way. To define the action of $\mathfrak S_n$ on $S\times_B B''$, we let $\mathfrak S_n$ act trivially on $S$ and act as permutations of components of $B''$; and to define the action on $(\mathbb P^1)^n\times B''$, we let $\mathfrak S_n$ act as permutations of components for both $(\mathbb P^1)^n$ and $B''$. It is clear from the construction that $\psi': S\times_BB''\dashrightarrow (\mathbb P^1)^n\times B''$ is $\mathfrak S_n$-equivariant. Thus $\psi'$ induces a rational map $ S\dashrightarrow ((\mathbb P^1)^n\times B'')/\mathfrak S_n$ over $B$. Let $S'$ be the image of this map. Thus we get a dominant rational map
\[\phi: S\dashrightarrow S'.
\]
\begin{prop}\label{PropOfPhi}
$S'$ is a surface and the degree of $\phi$ divides $d$. Furthermore, if $n\geq 2$, the general fiber of $S'\dashrightarrow B$ is of geometric genus $\leq (\frac{d}{\deg\phi}-1)^2$.
\end{prop}
\begin{prf}
Over the general point $b=(b_1,\ldots, b_n)\in B''$, $\psi'$ is given by the morphism $\psi'_b: C\to (\mathbb P^1)^n$ induced by the $n$ morphisms $C\to \mathbb P^1$ of degree $d$ corresponding to the points $b_i\in B'$, where $C$ is the fiber of $f: S\dashrightarrow B$ oveer $\pi'(b)\in B$. Let $C'$ be the image of $\psi_b'$. Then the fiber of $S'\dashrightarrow B$ over $\pi'(b)\in B$ is $C'$ by construction. Thus $S'$ is a surface, and the degree of $\phi$ is the degree of $C$ over $C'$, which divides $d$. This proves the first statement. To prove the second, we need to prove the geometric genus of $C'$ is $\leq (\frac{d}{\deg\phi}-1)^2$. Since $C'$ is a curve of degree $(\frac{d}{\deg\phi},\ldots, \frac{d}{\deg\phi})$ in $(\mathbb P^1)^n$, we can use the Castelnuovo type lemma below concerning algebraic curves in $(\mathbb P^1)^n$.
\end{prf}
\begin{lmm}\label{GenusOfCurvesInP1n}
Let $n\geq 2$. Let $C$ be an integral curve in $(\mathbb P^1)^n$ of degree $(d,\ldots, d)$. Then the geometric genus of $C$ is less than or equal to $(d-1)^2$.
\end{lmm}
\begin{prf}
We prove by induction on $n$. When $n=2$, it is the adjunction formula. Now assume that any curve $C'\subset (\mathbb P^1)^{n-1}$ of degree $(e,\ldots, e)$ has geometric genus $\leq (e-1)^2$. Consider the projection $C\to C''\subset (\mathbb P^1)^{n-1}$ to the first $n-1$ components. The degree of $C\to C''$ is of the form $d/e$, for some $e$. Hence,  $C''\subset (\mathbb P^1)^{n-1}$ is a curve of degree $(e,\ldots, e)$, hence it has geometric genus $\leq (e-1)^2$ by induction assumption. Let $\tilde C$ and $\tilde C''$ be the normalization of $C$ and $C''$ respectively. Then $\tilde C$ is birational to its image $C'''$ in $\tilde C''\times \mathbb P^1$, where the map to the second component is given by the composition map $\tilde C\to C\stackrel{i_n}{\to}\mathbb P^1$. Here, the map $i_n: C\to \mathbb P^1$ is the projection map to the $n$-th component. Note that $C'''$  is of degree $(d/e, d)$ in $\tilde C''\times \mathbb P^1$ and note that $NS(C''\times \mathbb P^1)=NS(C'')\oplus NS(\mathbb P^1)$. Adjunction formula and $g(C'')\leq (e-1)^2$ give us $p_a(C''')=\frac{d(g(C'')+d-1)}{e}-d+1\leq d(e+d/e-2)-d+1\leq (d-1)^2$ since $1\leq e\leq d$. Hence, $p_g(C)=g(\tilde C)=p_g(C''')\leq p_a(C''')\leq (d-1)^2$, as desired.
\end{prf}

\subsection{Proof of Theorem~\ref{ThmFibgonDeK3}; The Case (a)}\label{ProofCaseA}
Let $f: S\dashrightarrow B=\mathbb P^1$ be a fibration into curves realizing the fibering gonality of $S$. After a desingularization of indeterminacies of $f$, we get a dominant morphism $\tilde f: \tilde S\to B$ whose general fiber is of gonality $d=\mathrm{Fibgon}(S)$. In Section~\ref{ConstructionofPhi}, we constructed a surface $S'$ that is a fibration over $B$ into curves and a dominant rational map $\phi: \tilde S\dashrightarrow S'$ over $B$ of degree dividing $d$. The Kodaira dimension $\kappa(S')$, the irregularity $q(S')$ and the geometric genus $p_g(S')$ of $S'$ cannot exceed those of $S$ since $S'$ is dominated by $S$. By Enriques-Kodaira classification of algebraic surfaces~\cite{EnriquesKodaira}, $S'$ can only be birational to $\mathbb P^2$, a K3 surface or an Enriques surface. 

If  $S'$ is a rational surface, then $\mathrm{Irr}(S)\leq 
\deg \phi \leq d=\mathrm{Fibgon}(S)$. Here, the inequality $\deg\phi\leq d$ is because of Proposition~\ref{PropOfPhi}. But clearly $\mathrm{Fibgon}(S)\leq \mathrm{Irr}(S)$. We get the equality. This is case (a) of the theorem. We will treat the case when $S'$ is birational to a $K3$ surface in Section~\ref{ProofCaseB} and we will prove in Proposition~\ref{InequalityAsb} the inequality given in case (b) of the theorem. 

It remains to eliminate the case when $S'$ is birational to an Enriques surface $S''$. After a birational modification of $\tilde S$, there is a dominant morphism $g:\tilde S\to S''$. Since $\tilde S$ is simply-connected, $g$ factors through the universal covering $S'''$ of $S''$. $S'''$ is a K3 surface of Picard number at least $10$. But Lemma~\ref{b2trTheSame} below shows that the Picard number of $S'''$ can only be $1$ since it is dominated by $S$. This gives us a contradiction. The proof of Theorem~\ref{ThmFibgonDeK3} is thus concluded when $S'$ is either rational or an Enriques surface.

\begin{lmm}\label{b2trTheSame}
Let $\phi: X\dashrightarrow X'$ be a dominant rational map between projective hyper-Kähler manifolds of the same dimension. Then $b_{2,tr}(X)=b_{2,tr}(X')$.
\end{lmm}
\begin{prf}
Let $\tau: \tilde X\to X$ and $\tilde\phi:\tilde X\to X'$ be a resolution of indeterminacy points of $\phi$. Then $\tau^*: H^2(X,\mathbb Q)_{tr}\to H^2(\tilde X,\mathbb Q)_{tr}$ is an isomorphism and $\tilde\phi^*:H^2( X',\mathbb Q)_{tr}\to H^2(\tilde X,\mathbb Q)_{tr}$ is injective. They are moreover both morphisms of Hodge structures. But $H^2(X,\mathbb Q)_{tr}$ and $H^2(X',\mathbb Q)_{tr}$ are simple Hodge structures. This implies that $\tilde\phi^*:H^2( X',\mathbb Q)_{tr}\to H^2(\tilde X,\mathbb Q)_{tr}$ is an isomorphism and hence the result.
\end{prf}

The proof of Theorem~\ref{ThmFibgonDeK3} in the case where $S'$ is birational to a K3 surface will be completed in Section~\ref{ProofCaseB}.

\subsection{Rational maps between $K3$ surfaces}\label{RatMapsBetweenK3}
We treat rational maps between $K3$ surfaces in this part. Let $S$ (resp. $S'$) be a projective $K3$ surface whose Picard group is generated by an ample line bundle of degree $2D-2$ (resp. $2D'-2$). Let $\phi: S\dashrightarrow S'$ be a dominant rational map. 
\begin{prop}\label{BoundsForDegreeOfK3}
We have the folllowing inequality \[\frac{1}{(\deg\phi)^{21}}\leq \frac{D-1}{D'-1}\leq (\deg\phi)^{21}.\]
\end{prop}
\begin{prf}
Let $\tau: \tilde S\to S$ and $\tilde\phi: \tilde S\to S'$ be a resolution of indeterminacy points of $\phi: S\dashrightarrow S'.$ Let $T$ (resp. $T'$) be the lattice $H^2(S,\mathbb Z)_{tr}$ (resp. $H^2(S',\mathbb Z)_{tr}$) endowed with the intersection form. For a positive integer $e$, define $T'(e)$ to be the lattice $T'$ with the quadratic form multiplied by $e$. For example, with this notation, the sublattice $eT'$ of $T'$ is isometric to $T'(e^2)$ as lattices. The image $E$ of the morphism $\tilde\phi^*: H^2(S',\mathbb Z)_{tr}\to H^2(\tilde S,\mathbb Z)_{tr}\cong T$, viewed as a sublattice of $T$, is isometric to $T'(\deg\phi)$. The isomorphism $H^2(\tilde S,\mathbb Z)_{tr}\cong T$ is because  $\tau^*:H^2(S,\mathbb Z)_{tr}\to H^2(\tilde S,\mathbb Z)_{tr}$ is an isomorphism. We thus get the following equalities
\begin{equation}\label{EquationOfDisc}
    [T: E]^2=\left|\frac{\mathrm{disc}(E)}{\mathrm{disc}(T)}\right|=(\deg\phi)^{21}\left|\frac{\mathrm{disc}(T')}{\mathrm{disc}(T)}\right|=(\deg\phi)^{21}\frac{D'-1}{D-1}.
\end{equation}
On the other hand, we have
\begin{lmm}\label{philowerstarisinjective}
    The morphism of abelian groups $\tilde\phi_*: H^2(\tilde S,\mathbb Z)_{tr}\to H^2(S',\mathbb Z)_{tr}$ is injective and sends $E$ onto $(\deg\phi)T'$.
\end{lmm}
\begin{prf}
The projection formula shows that $\tilde\phi_*\tilde\phi^*=\deg\phi\cdot Id$. Hence, $\tilde\phi_*$ sends $E$ onto $(\deg\phi)T'$. By Lemma~\ref{b2trTheSame} and the fact that $\tilde\phi_*$ is surjective with $\mathbb Q$-coefficients, the kernel of $\tilde\phi_*$ is of torsion. But $H^2(\tilde S,\mathbb Z)_{tr}$ is torsion-free, as $\tilde S$ is simply connected. We conclude that the kernel of $\tilde\phi_*$ is zero, as desired. 
\end{prf}

Now Lemma~\ref{philowerstarisinjective} implies that  \begin{equation}\label{InequalityOfTE}
    1\leq [T:E]\leq [T':(\deg\phi)T']=(\deg\phi)^{21}.
\end{equation} 
Proposition~\ref{BoundsForDegreeOfK3} follows by combining (\ref{EquationOfDisc}) and (\ref{InequalityOfTE}).
\end{prf}
\begin{rmk}
One can similarly prove the more general result on hyper-Kähler manifolds, namely Proposition~\ref{InequalityOfLatticesHK}. The detailed proof is given in Section~\ref{Inequalities}. A sharper lower bound will also be given there.
\end{rmk}

\subsection{Proof of Theorem~\ref{ThmFibgonDeK3}; The Case (b)}\label{ProofCaseB}
Now let us continue the proof of Theorem~\ref{ThmFibgonDeK3}. Let $S$ be a $K3$ surface of Picard number $1$ and let $f: S\dashrightarrow B=\mathbb P^1$ be a fibration into curves realizing the fibering gonality. Let $\phi: S\dashrightarrow S'$ be the rational map constructed as in Section~\ref{ConstructionofPhi}. Recall that  $S'$ can only be birational to $\mathbb P^2$, a K3 surface or an Enriques surface, and we have treated, in Section~\ref{ProofCaseA}, the cases when $S'$ is birational to $\mathbb P^2$ or to an Enriques surface. In the rest of this section, we discuss the case when $S'$ is birational to a $K3$ surface. By changing the birational model, we may assume $S'$ \emph{is} a K3 surface. By Lemma~\ref{b2trTheSame}, the Picard number of $S'$ is also $1$. The following proposition shows that in our situation the Case (b) of Theorem~\ref{ThmFibgonDeK3} holds, which concludes the proof of Theorem~\ref{ThmFibgonDeK3}.
\begin{prop}\label{InequalityAsb}
The following inequality holds: $$\mathrm{Fibgen}(S)\leq \mathrm{Fibgon}(S)^{21/2}.$$
\end{prop}
\begin{prf}
Let $D$ and $D'$ be the degrees of the $K3$ surfaces $S$ and $S'$, respectively. Let $C'$ be the general fiber of $S'\dashrightarrow B$ as in Section~\ref{ConstructionofPhi}. Then we have the following inequalities
\begin{align*}
    \left(\frac{\mathrm{Fibgon}(S)}{\deg\phi}-1\right)^2 & \geq p_g(C')  & \textrm{ by Proposition~\ref{PropOfPhi}}\\
    & \geq \sqrt{\frac{D'}{2}}  & \textrm{ by Ein--Lazarsfeld's theorem (Theorem~\ref{EinLazarsfeld})}\\
    & \geq \sqrt{\frac{D}{2(\deg\phi)^{21}}}  & \textrm{ by Proposition~\ref{BoundsForDegreeOfK3}}\\
    & \geq \frac{\mathrm{Fibgen}(S)}{4(\deg\phi)^{21/2}}  & \textrm{ by Ein--Lazarsfeld's theorem (Theorem~\ref{EinLazarsfeld})}.
\end{align*}
Note that $\mathrm{Fibgon}(S)\geq 2\deg\phi$. Proposition~\ref{InequalityAsb} follows from these inequalities.
\end{prf}

\section{Some inequalities about Picard lattices of hyper-Kähler manifolds}\label{Inequalities}
We prove in this section Propositions~\ref{InequalityOfLatticesHK} and \ref{InequalityOfLatticesK3}. 

\begin{proof}[Proof of Proposition~\ref{InequalityOfLatticesHK}]
Let $\tau: \tilde X\to X$ and $\tilde\phi: \tilde X\to X'$ be a resolution of indeterminacy points of $\phi: X\dashrightarrow X'$. Let $T$ (resp. $T'$) be the lattice $H^2(X,\mathbb Z)_{tr}$ (resp. $H^2(X',\mathbb Z)_{tr}$) endowed with the Beauville--Bogomolov--Fujiki form. As in the proof of Proposition~\ref{BoundsForDegreeOfK3}, for a positive integer $e$, define $T'(e)$ to be the lattice $T'$ with the quadratic form multiplied by $e$. We claim that the image $E$ of the morphism $\tilde\phi^*: H^2(X',\mathbb Z)_{tr}\to \tau^*H^2(X,\mathbb Z)_{tr}\cong T$, viewed as a sublattice of $T$, is isometric to $T'((\deg\phi)^{\frac1n})$. This follows from the equalities $[q_X(\tilde\phi^*\alpha)]^n=c_X\cdot(\int_X\tilde\phi^*\alpha^{2n})=(\deg\phi)\cdot c_X\cdot(\int_{X'}\alpha^{2n})=(\deg\phi)\cdot [q_{X'}(\alpha)]^n$, where $c_X=c_{X'}$ is the Fujiki constant for the deformation class of $X$. Now the claim implies the following equalities
\begin{equation}\label{EquationOfDiscHK}
    [T: E]^2=\left|\frac{\mathrm{disc}(E)}{\mathrm{disc}(T)}\right|=(\deg\phi)^{\frac{b_{2,tr}(X)}{n}}\cdot\left|\frac{\mathrm{disc}(T')}{\mathrm{disc}(T)}\right|=(\deg\phi)^{\frac{b_{2,tr}(X)}{n}}\cdot\left|\frac{\mathrm{disc}(\mathrm{Pic}(X'))}{\mathrm{disc}(\mathrm{Pic}(X))}\right|.
\end{equation}
On the other hand, with a similar argument to Lemma~\ref{philowerstarisinjective}, we prove that the morphism of abelian groups $\tilde\phi_*: H^2(X,\mathbb Z)_{tr}\to H^2(X',\mathbb Z)_{tr}$ injective and sends $E$ to $(\deg\phi)T'$. Hence, \begin{equation}\label{InequalityOfTEHK}
    1\leq [T:E]\leq [T':(\deg\phi)T']=(\deg\phi)^{b_{2,tr}(X)}.
\end{equation} 
The proposition follows by combining (\ref{EquationOfDiscHK}) and (\ref{InequalityOfTEHK}).
\end{proof}

\begin{proof}[Proof of Proposition~\ref{InequalityOfLatticesK3}]
The only thing that needs proving, in the view of Proposition~\ref{InequalityOfLatticesHK}, is the following inequality
\[\frac{1}{(\deg\phi)^{\rho(X)}}\leq\left| \frac{\mathrm{disc}(\mathrm{Pic}(S))}{\mathrm{disc}(\mathrm{Pic}(S'))}\right|.
\]
Let $\tau: \tilde S\to S$ and $\tilde\phi: \tilde S\to S'$ be a resolution of indeterminacy points of $\phi: S\dashrightarrow S'$. The morphism $\tilde\phi^*: \mathrm{Pic}(S')\to \mathrm{Pic}(\tilde S)$ enlarges the quadratic form by $\deg\phi$ since $\tilde\phi^*\alpha\cup\tilde\phi^*\beta=\tilde\phi^*(\alpha\cup\beta)=\deg\phi\cdot(\alpha\cup\beta)$ for $\alpha, \beta\in \mathrm{Pic}(S')$. Thus 
\begin{equation}\label{EquationOfDiscK3enlargie}
    \mathrm{disc}(\tilde\phi^*(\mathrm{Pic}(S')))=(\deg\phi)^{\rho(S)}\mathrm{disc}(\mathrm{Pic}(S')).
\end{equation} 
\begin{lmm}\label{LatticesOrthogonalLemma}
The sublattice $\ker(\tilde\phi_*:\mathrm{Pic}(\tilde S)\to \mathrm{Pic}(S'))$ of $\mathrm{Pic}(\tilde S)$ is the orthogonal complement of $\tilde\phi^*(\mathrm{Pic}(S'))$ in $\mathrm{Pic}(\tilde S)$.
\end{lmm} 
\begin{prf}
Let $\alpha\in \mathrm{Pic}(\tilde S)$. Let us show that $\tilde\phi^*\alpha=0$ if and only if for any  $\beta\in\mathrm{Pic}(S')$, we have $\alpha\cup\tilde\phi^*\beta=0$ in $H^4(\tilde S, \mathbb Z)$. The projection formula gives $\tilde\phi_*(\alpha\cup\tilde\phi^*\beta)=(\tilde\phi_*\alpha)\cup\beta$ in $H^4(S', \mathbb Z)$. Hence, if $\tilde\phi^*\alpha=0$, then $\tilde\phi_*(\alpha\cup\tilde\phi^*\beta)=0$. But $\tilde\phi_*: H^4(\tilde S, \mathbb Z)\to H^4(S', \mathbb Z)$ is an isomorphism, we must have $\alpha\cup\tilde\phi^*\beta=0$. Conversely, if $\alpha\cup\tilde\phi^*\beta=0$ for every $\beta\in\mathrm{Pic}(S')$, still by the projection formula, we get $(\tilde\phi_*\alpha)\cup\beta=0$, which implies that $\tilde\phi_*\alpha=0$ since the intersection product map is nondegenerate on $\mathrm{Pic}(S')$.
\end{prf}

Taking into account the fact that the intersection map on $H^2(\tilde S,\mathbb Z)$ is nondegenerate on $\tilde\phi^*\mathrm{Pic}(S')$, Lemma~\ref{LatticesOrthogonalLemma} implies that $\ker(\tilde\phi_*)\oplus \tilde\phi^*(\mathrm{Pic}(S'))$ is a direct sum and that $\ker(\tilde\phi_*)\oplus \tilde\phi^*(\mathrm{Pic}(S'))$ is of \emph{finite} index in the Abelian group $\mathrm{Pic}(\tilde S)$.  Hence, 
\begin{equation}\label{InequalityOfTEK3}
    \left|\frac{\mathrm{disc}(\ker(\tilde\phi_*))\cdot\mathrm{disc}(\tilde\phi^*(\mathrm{Pic}(S'))}{\mathrm{disc(\mathrm{Pic}(\tilde S))}}\right|=[\mathrm{Pic}(\tilde S): \ker(\tilde\phi_*)\oplus \tilde\phi^*(\mathrm{Pic}(S'))]^2.
\end{equation}
Since $\tilde\phi_*:\mathrm{Pic}(\tilde S)\to \mathrm{Pic}(S'))$ sends $\ker(\tilde\phi_*)\oplus \tilde\phi^*(\mathrm{Pic}(S'))$ onto $(\deg\phi)\mathrm{Pic}(S')$, and since the induced morphism \[\overline{\tilde\phi_*}: \mathrm{Pic}(\tilde S)/\ker\tilde\phi_*\to \mathrm{Pic}(S')\] is injective, we have 
\begin{equation}\label{InequalityOfThingsK3}
    [\mathrm{Pic}(\tilde S): (\ker(\tilde\phi_*)\oplus \tilde\phi^*(\mathrm{Pic}(S'))]=[(\mathrm{Pic}(\tilde S)/\ker\tilde\phi_*): \overline{\tilde\phi^*(\mathrm{Pic}(S'))}]\leq [\mathrm{Pic}(S'):(\deg\phi)\mathrm{Pic}(S')]=(\deg\phi)^{\rho(S)},
\end{equation}
where $\overline{\tilde\phi^*(\mathrm{Pic}(S'))}$ is the image of $\tilde\phi^*(\mathrm{Pic}(S'))$ in $\mathrm{Pic}(\tilde S)/\ker\tilde\phi_*$.
Notice the following Lemma.
\begin{lmm}\label{BlowUpDoesntChangeDisc}
$|\mathrm{disc}(\mathrm{Pic}(\tilde S))|=|\mathrm{disc}(\mathrm{Pic}(S))|$
\end{lmm}
\begin{prf}
$\tilde S$ is obtained by a sequence of blowing-ups of points from $S$. Therefore, 
\[\mathrm{Pic}(\tilde S)=\tau^*\mathrm{Pic}(S)\oplus \bigoplus_i\mathbb Z E_i,\]
where $E_i$ is the \emph{total transform} in $\tilde S$ of the exceptional divisor of the $i$-th blowing-up. We have the following formula for the intersection numbers of $E_i$:
\[E_i.E_j=\left\{\begin{array}{cc}
   -1  & \textrm{ if } i=j \\
   0  & \textrm{ otherwise. }
\end{array}\right.
\]
Hence, $\mathrm{disc}(\mathrm{Pic}(\tilde S))=\mathrm{disc}(\mathrm{Pic}( S))\cdot\mathrm{disc}(\oplus_i\mathbb ZE_i)=\pm\mathrm{disc}(\mathrm{Pic}(S))$, as desired.
\end{prf}

The proposition now follows from Lemma~\ref{BlowUpDoesntChangeDisc} and from inequalities (\ref{EquationOfDiscK3enlargie}), (\ref{InequalityOfTEK3}) and (\ref{InequalityOfThingsK3}), noticing that $\mathrm{disc}(\ker(\tilde\phi_*))\geq 1$.
\end{proof}

 ~\newline
 
Sorbonne Université and Université de Paris, CNRS, IMJ-PRG, F-75005 Paris, France.

\textit{Email adress}: \verb|chenyu.bai@imj-prg.fr|.


\begin{thebibliography}{}
\bibitem{Arbarello} E. Arbarello, M. Cornalba, P. A. Griffiths, J. Harris, \emph{Geometry of Algebraic Curves, Volume I}, 1985, Springer. 

\bibitem{Irrationality} F. Bastianelli, P. De Poi, L. Ein, R. Lazarsfeld, B. Ullery, \emph{Measures of irrationality for hypersurfaces of large degree}, Compos. Math. 153 (2017), no. 11, 2368-2393.

\bibitem{EinLazarsfeld} L. Ein, R. Lazarsfeld, \emph{The Konno invariant of some algebraic varieties}, European Journal of Mathematics volume 6, 420–429 (2020).

\bibitem{EnriquesKodaira} K. Kodaira, (1964), \emph{On the structure of compact complex analytic surfaces. I-IV}, American Journal of Mathematics, 86 (4): 751–798, 88 (3): 682–721, 90 (1): 55–83, 90 (4): 1048–1066. 1964,1966, 1968.

\bibitem{ThesisOfMartin} O. Martin, \emph{Zero-cycles and measures of irrationality for abelian varieties}, PhD Thesis, (2020).

\bibitem{ThesisOfStapleton} D. Stapleton, \emph{The degree of irrationality of very general hypersurfaces in some homogeneous spaces}, Ph.D. thesis, Stony Brook University, 2017.

\bibitem{Voisin} C. Voisin. \emph{Théorie de Hodge et géométrie algébrique complexe}, Cours spécialisés, SMF 2003.

\bibitem{IrrAVVoisin} C. Voisin, \emph{Chow ring and the gonality of general abelian varieties}, Ann. H. Lebesgue 1 (2018), 313–332.

\bibitem{FibgenFibgon} C. Voisin, \emph{On fibrations and measures of irrationality of hyper-Kähler manifolds}, la  Revista de la Unión Matemática Argentina, volume spécial MCA 2021.



\bibitem{Yang} R. Yang, \emph{On irrationality of hypersurfaces in $\mathbb P^{n+1}$}, Proc. Amer. Math. Soc. 147 (2019), 971–976.

\end{thebibliography}
\end{document}